\pdfoutput=1
\documentclass[aoas]{imsart}
\usepackage{graphicx}
\RequirePackage[OT1]{fontenc}
\RequirePackage{amsthm,amsmath,natbib}
\RequirePackage[colorlinks]{hyperref}
\RequirePackage{hypernat}

\startlocaldefs
\newcommand{\Supp}{\operatorname{Supp}}
\newcommand{\Int}{\operatorname{Int}}
\newcommand{\beqast}{\begin{eqnarray*}}
\newcommand{\eeqast}{\end{eqnarray*}}

\numberwithin{equation}{section}
\theoremstyle{plain}

\newtheorem{theorem}{Theorem}[section]
\newtheorem{definition}{Definition}[section]
\newtheorem{claim}{Claim}[section]
\def\eref#1{(\ref{#1})}

\newcommand{\aml}{\underset{\ell\in L}{\arg\min}}
\newcommand{\amal}{\underset{\ell\in L}{\arg\max}}
\newcommand{\beq}{\begin{equation}}
\newcommand{\eeq}{\end{equation}}
\newcommand{\beqa}{\begin{eqnarray}}
\newcommand{\eeqa}{\end{eqnarray}}
\newcommand{\ba}{\begin{array}}
\newcommand{\ea}{\end{array}}
\newcommand{\bac}{\begin{array}{ccccccccccc}}
\newcommand{\eac}{\end{array}}
\newcommand{\co}[1]{}
\endlocaldefs

\begin{document}
\begin{frontmatter}
\title{A Principal Component Analysis for Trees}
\runtitle{Tree-Line Analysis}
\begin{aug}
\author{\fnms{Burcu} \snm{Ayd{\i}n}\thanksref{t1}\ead[label=e1]{aydin@email.unc.edu}},
\author{\fnms{G\'{a}bor} \snm{Pataki}\ead[label=e2]{pataki@email.unc.edu}},
\author{\fnms{Haonan} \snm{Wang}\thanksref{t3} \ead[label=e3]{wanghn@stat.colostate.edu}},
\author{\fnms{Elizabeth} \snm{Bullitt}\thanksref{t4} \ead[label=e4]{bullitt@med.unc.edu}}
\and
\author{\fnms{J.S.} \snm{Marron}\thanksref{t5} \ead[label=e5]{marron@email.unc.edu}}

\thankstext{t1}{Partially supported by NSF grant DMS-0606577, and NIH Grant RFA-ES-04-008.}
\thankstext{t3}{Partially supported by NSF grant DMS-0706761.}
\thankstext{t4}{Partially supported by NIH grants R01EB000219-NIH-NIBIB and R01 CA124608-NIH-NCI.}
\thankstext{t5}{Partially supported by NSF grant DMS-0606577, and NIH Grant RFA-ES-04-008.}

\affiliation{University of North Carolina and Colorado State University}

\address{Department of Statistics and Operations Research\\
University of North Carolina\\
Chapel Hill, NC 27599-3260\\
\printead{e1}\\
\phantom{E-mail:\ }\printead*{e2}\\
\phantom{E-mail:\ }\printead*{e5}}

\address{Department of Statistics\\
Colorado State University\\
Fort Collins, CO 80523-1877\\
\printead{e3}}

\address{Department of Surgery\\
University of North Carolina\\
Chapel Hill, NC 27599-3260\\
\printead{e4}}

\end{aug}

\begin{abstract}
The active field of Functional Data Analysis (about understanding the
variation in a set of curves) has been recently extended to Object Oriented
Data Analysis, which considers populations of more general objects.
A particularly challenging extension of this set of ideas is to populations of
tree-structured objects.
We develop an analog of Principal Component Analysis for trees, based on the
notion of tree-lines, and propose numerically fast (linear time) algorithms to solve
the resulting optimization problems.  The solutions we obtain are used in the analysis
of a data set of $73$ individuals, where each data object is a tree of blood vessels in
one person's brain.
\end{abstract}


\end{frontmatter}

\section{Introduction} \label{Sec:Intro}

Functional data analysis has been a recent active research area. See Ramsay
and Silverman (2002, 2005) for a good introduction and overview, and Ferraty and
Vieu (2006) for a more recent viewpoint. A major difference between this
approach, and more classical statistical methods is that curves are viewed as
the \textit{atoms} of the analysis, i.e. the goal is the statistical analysis
of a \textit{population of curves}.

Wang and Marron (2007) recently extended functional data analysis to Object
Oriented Data Analysis (OODA), where the atoms of the analysis are allowed to
be more general data objects. Examples studied there include images, shapes
and tree structures as the atoms, i.e. the basic data elements of the
population of interest. Other recent examples are populations of movies,
such as are being subjects of functional magnetic resonance imaging. A major
contribution of Wang and Marron (2007) was the development of a set of
tree-population analogs of standard functional data analysis techniques, such
as Principal Component Analysis (PCA). The foundations were laid via the
formulation of particular optimization problems, whose solution resulted in
that analysis method (in the same spirit in which ordinary PCA can be
formulated in terms of optimization problems).

Here the focus is on the challenging OODA case of tree structured data
objects. A limitation of the work of Wang and Marron (2007) was that no
general solutions appeared to be available for the optimization problems that
were developed. Hence, only limited toy examples (three and four node trees,
which thus allowed manual solutions) were used to illustrate the main ideas
(although one interesting real data lesson was discovered even with that
strong limitation on tree size).

One of our main contributions is that, through a detailed analysis of
the underlying optimization problem, and a complete solution of it, a linear
time computational method is now available. This allows the first actual
OODA of a production scale data set of a population of tree structured
objects. Ideas are illustrated in Section \ref{Sec:DatAnal} using a set of blood vessel trees in the
human brain, collected as described in Aylward and Bullitt (2002). In
the present paper, we choose to consider only variation in the
\textit{topology} of the trees, i.e. we consider only the branching structure
and ignore other aspects of the data, such as location, thickness and
curvature of each branch.

Even with this topology only restriction, there is still an
important \textit{correspondence} decision that needs to be made: which branch should be put on the left, and which one
on the right, see Section \ref{Sec:Corr}. Later analysis will also include
location, orientation and thickness information, by adding attributes to the
tree nodes being studied. A useful set of ideas for pursuing that type of
analysis was developed by Wang and Marron (2007).

In Subsection \ref{Sec:Lines} we define our main data analytic concept, the \emph{tree-line}, and the notion of principal components
based on tree-lines. Here we also state, and illustrate our main result, Theorem \ref{main}, which will allow us to quickly compute the principal components.
Subsection \ref{Sec:Results} is devoted to our data analysis  using the blood vessel data: we
carefully compare the correspondence approaches, and present our findings based on the computed principal components.
In Section \ref{Sec:Opt} we prove Theorem \ref{main} along with a host of necessary claims.

\section{Data and Analysis}\label{Sec:DatAnal}

\begin{figure}
[ptb]
\begin{center}
\includegraphics[
natheight=371.437500pt,
natwidth=325.875000pt,
height=224.875pt,
width=197.5625pt
]%
{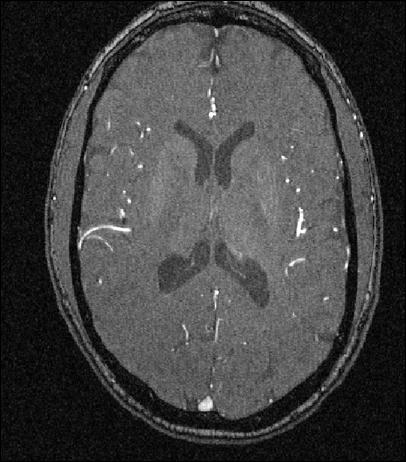}
\caption{Single Slice from a Magnetic Resonance Angiography image for one
patient. \ Bright regions indicate blood flow.}
\label{Fig1}%
\end{center}
\end{figure}

The data analyzed here are from a study of Magnetic Resonance Angiography brain images of a set
of $73$ human subjects of both sexes, ranging in age from $18$ to $72$, which can be found at Handle (2008).
One slice of one such image is shown in Figure \ref{Fig1}. This mode of imaging indicates
strong blood flow as white. These white regions are tracked in $3$ dimensions,
then combined, to give trees of brain arteries.

\begin{figure}
[ptb]
\begin{center}
\includegraphics[
natheight=231.250000pt,
natwidth=330.750000pt,
height=163.875pt,
width=233.5pt
]%
{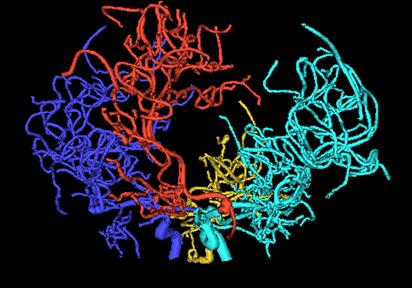}%
\caption{Reconstructed set of trees of brain arteries for the same patient
as shown in Figure \ref{Fig1}. The colors indicate regions of the brain: Gold (back),
Right (blue), Front (red), Left (cyan).}
\label{Fig2}%
\end{center}
\end{figure}

The set of trees developed from the image of which Figure \ref{Fig1} is one slice is
shown in Figure \ref{Fig2}. Trees are colored according to region of the brain. Each region is studied separately,
where each tree is one data point in the data set of its region. The goal of the present OODA
is to understand the population structure of $73$ subjects through $3$ data sets extracted from them: Back data set (gold trees), left data set (cyan) and right data set (blue).  
One point to note is that the front trees (red) are not studied here. This is because the source of flow for the front trees is variable, therefore this subpopulation has less biological meaning. For simplicity we chose to omit this sub-population.

The stored information for each of these trees is quite rich (enabling the
detailed view shown in Figure \ref{Fig2}). Each colored tree consists of a set of
branch segments. Each branch segment consists of a sequence of spheres fit
to the white regions in the MRA image (of which Figure \ref{Fig1} was one slice), as
described in Aylward and Bullitt (2002). Each sphere has a center (with $x$,
$y$, $z$ coordinates, indicating location of a point on the center line of the
artery), and a radius (indicating arterial thickness).

\subsection{Tree Correspondence\label{Sec:Corr}}

\begin{figure}
[ptb]
\begin{center}
\includegraphics[
natheight=408.500000pt,
natwidth=395.125000pt,
height=206.25pt,
width=199.5625pt
]%
{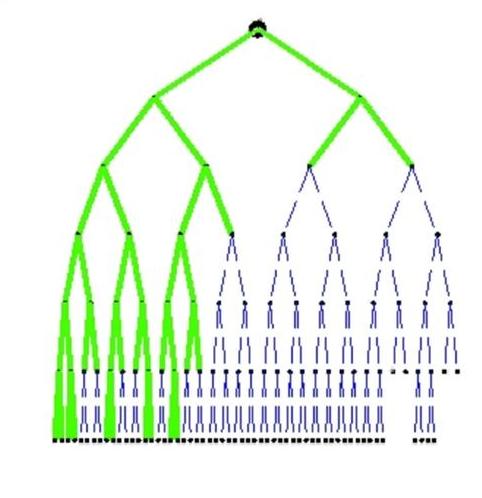}%
\caption{Green line segments show the topology only representation
of the gold (back tree) from Figure \ref{Fig2}. Only branching information is
retained for the OODA. Branch location and thickness information is
deliberately ignored. Thin blue curve shows the union over all trees in the
sample.}
\label{Fig3}
\end{center}
\end{figure}

Given a single tree, for example the gold colored (back) tree in Figure
\ref{Fig2}, we reduce it to only its topological (connectivity) aspects by
representing it as a simple binary tree. Figure \ref{Fig3} is an example of
such a representation. Each node in Figure \ref{Fig3} is best thought of as
a branch of the tree, and the green line segments simply show which child
branch connects to which parent. The root node at the top represents the
initial fat gold tree trunk shown near the bottom of Figure \ref{Fig2}. The
thin blue lines show the support tree, which is just the union of all of the
back trees, over the whole data set of $73$ patients.

There is one set of ambiguities in the construction of the binary tree shown
in Figure \ref{Fig3}. That is the choice, made for each adult branch, of which
child branch is put on the left, and which is put on the right. The
following two ways of resolving this ambiguity are considered here. Using
standard terminology from image analysis, we use the word
\textit{correspondence} to refer to this choice.

\begin{itemize}
\item \textbf{Thickness Correspondence:} Put the node that corresponds to
the child with larger median radius (of the sequence of spheres fit to the MRA
image) on the left. Since it is expected that the fatter child vessel will
transport the most blood, this should be a reasonable notion of
\textit{dominant branch}.

\item \textbf{Descendant Correspondence:} Put the node that corresponds to
the child with the most descendants on the left.
\end{itemize}

These correspondences are compared in Subsection \ref{Sec:Results}.

Other types of correspondence, that have not yet been studied, are also
possible. An attractive possibility, suggested in personal discussion by
Marc Niethammer, is to use location information of the children in this
choice. E.g. in the back tree, one could choose the child which is
physically more on the left side (or perhaps the child whose descendants are
more on average to the left) as the left node in this representation. This
would give a representation that is physically closer to the actual data,
which may be more natural for addressing certain types of anatomical issues.

\subsection{Tree-Lines\label{Sec:Lines}}

In this section we develop the tools of our main analysis, based on the
notion of \textit{tree-lines}. We follow the ideas of Wang and Marron
(2007), who laid the foundations for this type of analysis, with a set of
ideas for extending the Euclidean workhorse method of PCA to data sets of tree
structured objects. The key idea (originally suggested in personal
conversation by J. O. Ramsay) was to define an appropriate \textit{one
dimensional representation}, and then find the one
that best fits the data. The tree-line is a first simple approach to this problem.

First we define a binary tree:

\begin{definition}
A \textbf{binary tree} is a set of nodes that are connected by edges in a
directed fashion, which starts with one node designated as \textbf{root},
where each node has at most two children.
\end{definition}

\begin{figure}
[ptb]
\begin{center}
\includegraphics[
natheight=207.062500pt,
natwidth=484.812500pt,
height=64.125pt,
width=250pt
]
{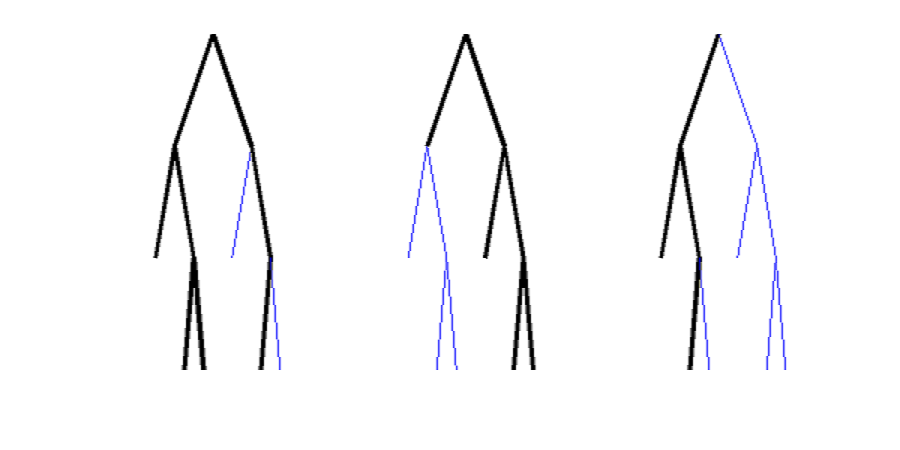}
\\
\caption{Toy example of a data set of trees, $T$, with $n=3$. This will
be used to illustrate several issues below.}
\label{Fig4}
\end{center}
\end{figure}

Using the notation $t_{i}$ for a single tree, we let
\beq \label{Tti}
T=\left\{  t_{1}%
,...,t_{n}\right\}
\eeq
denote a data set of $n$ such trees.
A toy example of a set of $3$ trees is given in Figure \ref{Fig4}.

To identify the nodes within each tree more easily, we use the level-order indexing method from Wang and Marron (2007).
The root node has index $1$. For the remaining nodes, if a node has index $\omega$, then the index of its left child
is $2 \omega$ and of its right child is $2 \omega+1$.
These indices enable us to identify a binary tree by only listing the indices of its nodes.

The basis of our analysis is an appropriate metric, i.e. distance, on tree
space. We use the common notion of Hamming distance for this purpose:

\begin{definition}
Given two trees $t_{1}$ and $t_{2}$, their  \textbf{distance} is
\[
d\left(  t_{1},t_{2}\right)  = | t_{1}\backslash t_{2} | + |  t_{2}\backslash t_{1} |,
\]
where $\backslash$ denotes set difference.
\end{definition}

Two more basic  concepts are defined below; the notion of support tree has already been shown in Figure \ref{Fig3} (as
the thin blue lines).

\begin{definition}
For a data set $T$, given as in \eref{Tti}, the support tree, and the intersection tree are defined as
\beqast
\Supp(T) & = & \cup_{i=1}^n t_i \\
\Int(T) & = & \cap_{i=1}^n t_i.
\eeqast
\end{definition}

Figure \ref{Fig6} shows the support trees of the data sets used in this study. Figure \ref{Fig7} includes the corresponding intersection trees.

The main idea of a tree-line (our notion of one dimensional representation) is
that it is constructed by adding a sequence of single nodes, where each new
node is a child of the most recent child:

\begin{definition}
A\textbf{\ tree-line}, $L=\left\{  \ell_{0},\cdots,\ell_{m}\right\}  $, is a sequence
of trees where $\ell_0$ is called the starting tree, and $\ell_{i}$ comes from $\ell_{i-1}$ by the addition of a single
node, labeled $v_{i}$. \ In addition each $v_{i+1}$ is a child of $v_{i}$.
\end{definition}

\begin{figure}
[ptb]
\begin{center}
\includegraphics[
natheight=207.062500pt,
natwidth=484.812500pt,
height=64.125pt,
width=300pt
]
{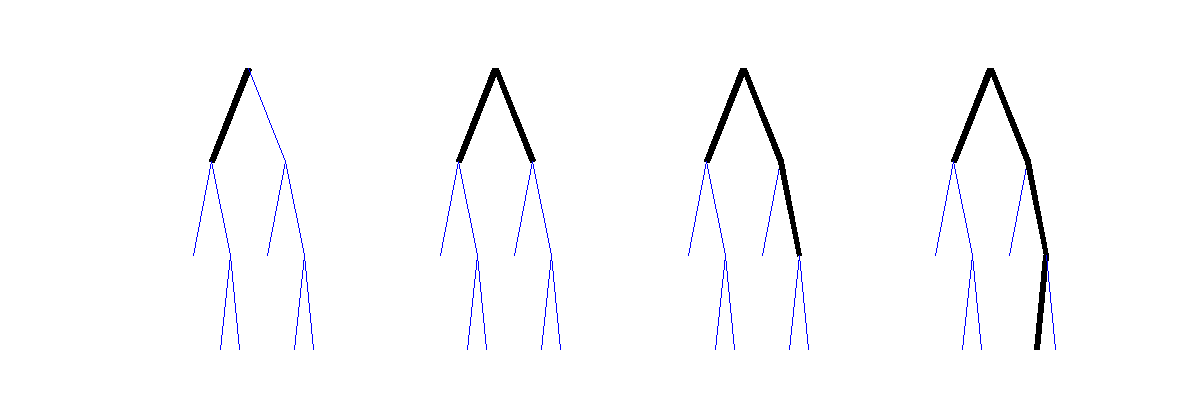}
\caption{Toy example of a tree-line. Each member come from
adding a node to the previous. Each new node is a child of the previously
added node. Starting point ($\ell_0$) is the intersection tree of the toy data set of Figure \ref{Fig4}.}
\label{Fig5}
\end{center}
\end{figure}

An example of a tree-line is given in Figure \ref{Fig5}.
Insight as to how well a given tree-line fits a data set is based upon the
concept of projection:

\begin{definition}
Given a data tree $t$, its \textbf{projection} onto the tree-line $L$ is%
\[
P_{L}\left(  t\right)  =\aml \{ d\left(t,\ell \right)\}  .
\]

\end{definition}

Wang and Marron (2007) show that this projection is always unique. This will also follow from Claim \ref{CLfirst} in Section \ref{Sec:Opt}, whose
characterization of  the projection will be the key in computing the principal component tree-lines, defined shortly.

The above toy examples provide an illustration. Let $t_{3}$ be
the third tree shown in Figure \ref{Fig4}. Name the trees in the
tree-line, $L$, shown in Figure \ref{Fig5}, as $\ell_0$,$\ell_1$,$\ell_2$,$\ell_3$. \ The set of distances
from $t_{3}$ to each each tree in $L$ is tabulated as%
\[
\begin{tabular}
[c]{lllll}%
$j$ & $0$ & $1$ & $2$ & $3$\\
$d\left(  t_{3},\ell_{j}\right)  $ & $6$ & $5$ & $4$ & $5$%
\end{tabular}
\]

The minimum distance is $4$, achieved at $j=2$, so the projection of $t_{3}$
onto the tree-line $L$ is $\ell_{2}$.

Next we develop an analog of the first principal component ($PC1$), by finding the
tree-line that best fits the data.
\begin{definition} \label{def-l1star}
For a data set $T$, the \textbf{first principal component tree-line}, i.e. $PC1$, is
\[
L_1^* \, = \, \underset{L}{\arg\min} \sum_{t_i \in T} d(t_i,P_L(t_i))
\]
\end{definition}
In conventional Euclidean PCA, additional components are restricted to lie in
the subspace orthogonal to existing components, and subject to that
restriction, to fit the data as well as possible. For an analogous notion in tree space, we first
need to define the concept of the union of tree-lines, and of a projection onto it.
\begin{definition}
Given tree-lines $L_1 = \{ \ell_{1,0}, \ell_{1,1}, \dots, \ell_{1,p_1} \}, \dots, L_q = \{ \ell_{q,0}, \ell_{q,1}, \dots, \ell_{q, p_q} \}$, their
\textbf{union} is the set of all possible unions of members of $L_1$ through $L_q$:
\begin{eqnarray*}
L_1 \cup \dots \cup L_q & = & \{ \ell_{1, i_1} \cup \dots \cup \ell_{q, i_q} \, | \, i_1 \in \{ 0, \dots, p_1 \}, \dots, i_q \in \{ 0, \dots, p_q \}. \}
\end{eqnarray*}
Given a data tree $t$, the projection of $t$ onto $L_1 \cup \dots \cup L_q$ is
\beq \label{projunion}
P_{L_1 \cup \dots \cup L_q}\left(  t \right)  =\underset{\ell\in L_1 \cup \dots \cup L_q}{\arg\min}\{ d\left(
t,\ell \right)\}.
\eeq
\end{definition}
In our non-Euclidean tree
space, there is no notion of orthogonality available, so we instead just ask
that the $2$nd tree-line fit as much of data as possible, when used in
combination with the first, and so on.
\begin{definition} \label{k}
For $k \geq 1$ the $k$th principal component tree-line is defined recursively as
\beq \label{lkstar}
L_k^*= \aml \sum_{t_i \in T} d(t_i,P_{L_1^*\cup \cdots\cup L_{k-1}^*\cup L}(t_i)),
\eeq
and it is abbreviated as $PCk$.
\end{definition}

For the concept of PC tree-lines to be useful, it is of crucial importance to be able to compute them efficiently. We need
another notion.

\begin{figure}
[ptb]
\begin{center}
\includegraphics[
natheight=207.062500pt,
natwidth=484.812500pt,
height=150pt,
width=200pt
]%
{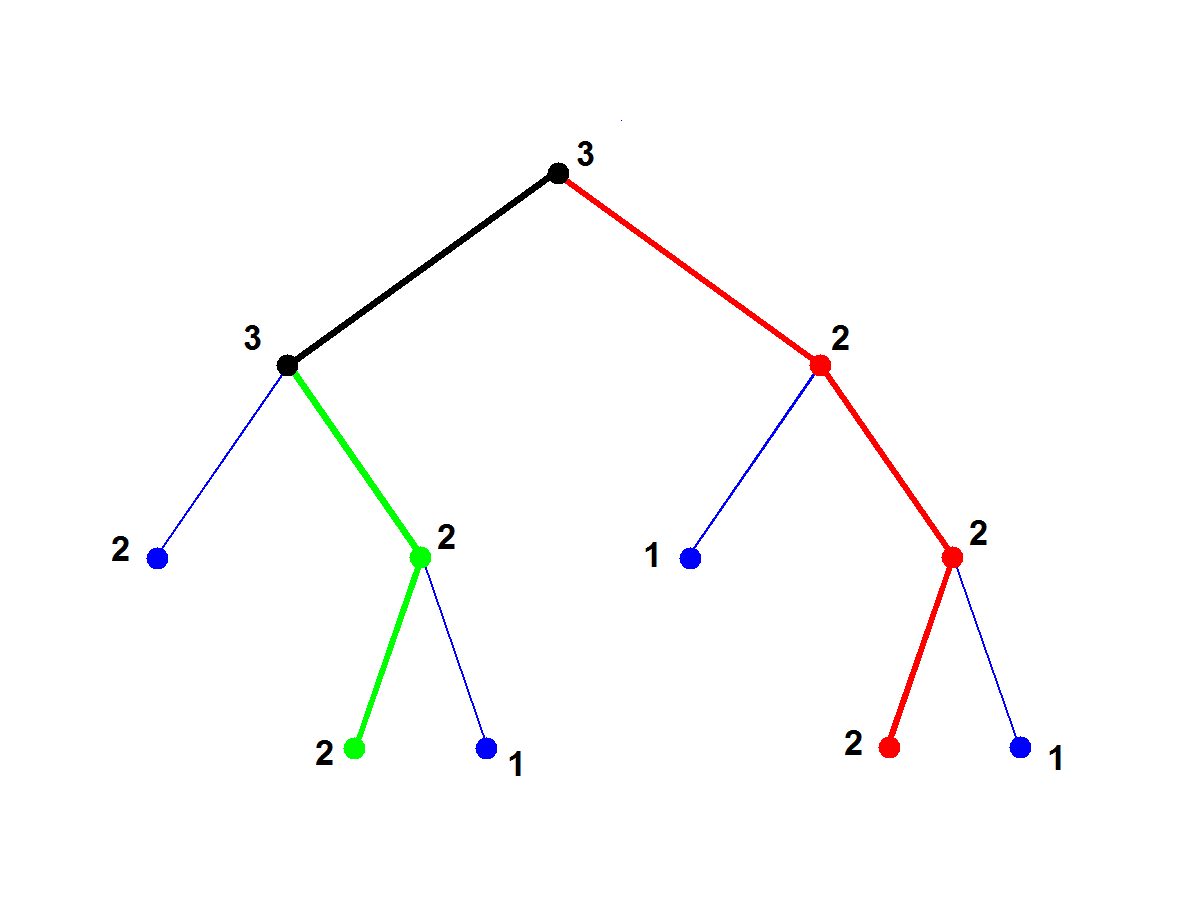}
\caption{Weighted support tree illustrating Theorem \ref{main}}
\label{Fig-wtedsupporttree}%
\end{center}
\end{figure}

\begin{definition}
Given a tree-line
$$
L=\left\{  \ell_{0}, \ell_1, \cdots, \ell_m \right\}
$$
we define the {\em path of} $L$ as
$$
V_L = \ell_m \setminus \ell_0.
$$
\end{definition}
Intuitively, a tree-line that well fits the data ``should grow in the direction that captures the most information''.
Furthermore, the $k$th PC tree-line should only aim to capture information that has not been explained by the first $k-1$ PC tree-lines.
This intuition is made precise in the following theorem, which is the main theoretical result of the paper:

\begin{theorem} \label{main}
Let $ k \geq 1, \,$ and $L_1^*, \dots, L_{k-1}^*$ be the first $k-1$ PC tree-lines.
For $v \in \Supp(T)$ define
\beq
\ba{rcl}
w_k(v) & = & \left\{ \ba{rl}  0, & \text{if} \,\, v \in  V_{L_1^*} \cup \dots \cup V_{L_{k-1}^*}, \\
                              \sum_{v \in t_i} 1, & \text{otherwise}
                      \ea \right.
\ea
\eeq
Then the $k$th PC tree-line $L_k^*$ is the tree-line whose path maximizes the sum of $w_k$ weights in the support tree, i.e.
$\sum_{v \in V_{L_k}^*} w_k(v)$.
\co{\beq
L_k \sum_{v \in V_L} w(v,T,S)
\eeq}
\end{theorem}
The proof of Theorem \ref{main} is given in Section \ref{Sec:Opt}. Figure \ref{Fig-wtedsupporttree} is an illustration: the weight of a node is
the number of times the node appears in the trees of Figure  \ref{Fig4}. The black edge is the intersection tree of the same data set. The maximum weight path
attached to $\Int(T)$ is the red path, which gives rise to the tree-line of Figure \ref{Fig5}, which is thus the first principal component of the data set of Figure \ref{Fig4}.

After setting the weights of the nodes on the red path to zero, the maximum weight path attached to Int(T) becomes the green path, which by Theorem \ref{main} gives rise to $PC2$. The usefulness of these tools is demonstrated with actual
data analysis of the full tree data set.

\subsection{Real Data Results\label{Sec:Results}}

This section describes an exploratory data analysis of the set of $n=73$ brain
trees discussed above using these tree-line ideas. The principal component
tree-lines are computed as defined in Theorem \ref{main}. Both
correspondence types, defined in Section \ref{Sec:Corr} are considered and
compared.

The different brain location types (shown as different colors in Figure
\ref{Fig2}) are analyzed as separate populations (i.e. the $n=73$ blue trees
are first considered to be the population, then the $n=73$ gold trees, etc.),
called \textit{brain location sub-populations}. This reveals some
interesting contrasts between the brain location types in terms of symmetry.

\begin{figure}
[h]
\begin{center}
\includegraphics[
natheight=9.375400in,
natwidth=12.500000in,
height=4in,
width=5in
]%
{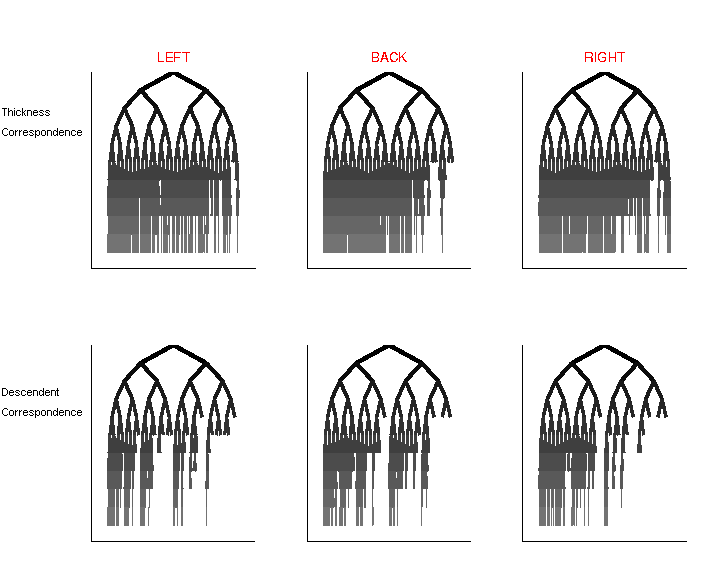}
\caption{Support trees, for both types of correspondence (shown in the
rows), and for three brain location tree types (shown in columns,
corresponding to the colors in Figure \ref{Fig2}). Shows that the descendant
correspondence gives a population with more compact variation than the
thickness correspondence.}
\label{Fig6}%
\end{center}
\end{figure}

We first compare the two types of correspondence defined in Section \ref{Sec:Corr}\ using the concept of the support tree. This
is done by displaying the support trees each type of correspondence, and
for each of the three tree location types (shown with different colors in Figure \ref{Fig2}), in
Figure \ref{Fig6}. Note that all of the support trees for the descendant
correspondence (bottom) are much smaller than for the thickness correspondence
(top), indicating that the descendant correspondence results in a much more
compact population. This seems likely to make it easier for our PCA method
to find an effective representation of the descendant based population.

Figure \ref{Fig6} already reveals an aspect of the population that was
previously unknown: there is not a very strong correlation
between median tree thickness of a branch, and the number of children.

\begin{figure}
[h]
\begin{center}
\includegraphics[
natheight=6.916800in,
natwidth=10.666600in,
height=3.7766in,
width=5.1232in
]%
{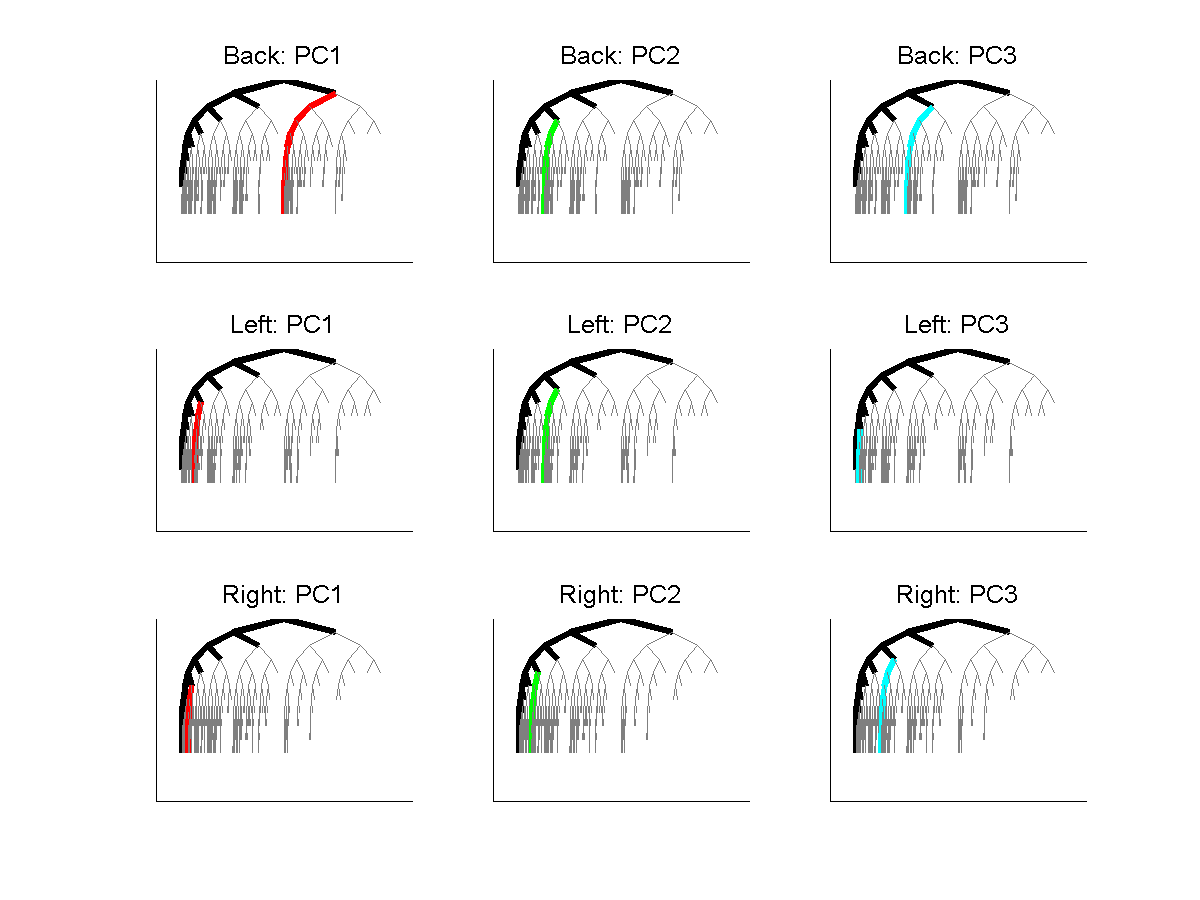}
\caption{Best fitting tree-lines, for different sub-populations (rows), and
PC number (columns). Intersection trees are shown in black.}
\label{Fig7}%
\end{center}
\end{figure}

Figure \ref{Fig7}\ shows the first $3$ PC tree-lines, for the three
sub-populations (shown as rows), with the intersection tree as the starting tree,
for the descendant correspondence.

In the human brain, the back circulation (gold) arises from a single vessel (the basilar artery) and immediately splits into two main trunks, supplying the back sides of the left and right hemispheres. These two parts of the back circulation are expected to be approximately mirror-image symmetrical with both sides containing one main vessel and other branches stemming from that. Consequently, for each tree on the back data set if we imagine a vertical axis that goes through the root node, we expect the subtrees on both sides of the axis to be symmetrical with each other.

The results of our model for the back subpopulation are consistent with this expectation. The main vessel of one of the hemispheres can be seen in the starting point (intersection tree) as the leftmost set of nodes, while the other main vessel becomes the first principal component.

As for the left and right circulations (cyan and blue trees) of the brain, they are expected to be close to mirror images of each other. Unlike the case of the back subpopulation, in each of these circulations there is a single trunk from which smaller branches stem. For this reason the bilateral symmetry observed within the back trees is not expected to be found here.

The fact that $PC1$'s for left and right subpopulations are at later splits suggest that the earlier splits tend to have relatively few descendants.
The remaining $PC2$ and $PC3$ tree-lines do not contain much additional information by themselves. However, when we consider PC's $1$,$2$ and $3$ together
and compare left and right subpopulations,i.e. compare the second and third rows of Figure \ref{Fig7}, the structural likeliness is quite visible.
It should also be noted that for both of the subpopulations all PC's are on the left side of the root-axis, indicating a strong bilateral asymmetry,
as expected.

The tree-lines, and insights obtained from them, were essentially similar for
the thickness correspondence, so those graphics are not shown here.

Next we study the tree-line analog of the familiar \textit{scores plot} from
conventional PCA (a commonly used high dimensional visualization device,
sometimes called a \textit{draftsman's plot} or a \textit{scatterplot
matrix}). In that case, the scores are the projection coefficients, which
indicate the size of the component of each data point in the given
eigen-direction. Pairwise scatterplots of these often give a set of useful
two dimensional views of the data. In the present case, given a data point
and a tree-line, the corresponding \textit{score} is just the length (i.e. the
number of nodes) of the projection. Unlike conventional PC scores, these are
all integer valued.

\begin{figure}
[h]
\begin{center}
\includegraphics[
natheight=6.916800in,
natwidth=10.666600in,
height=3.7766in,
width=5.1232in]
{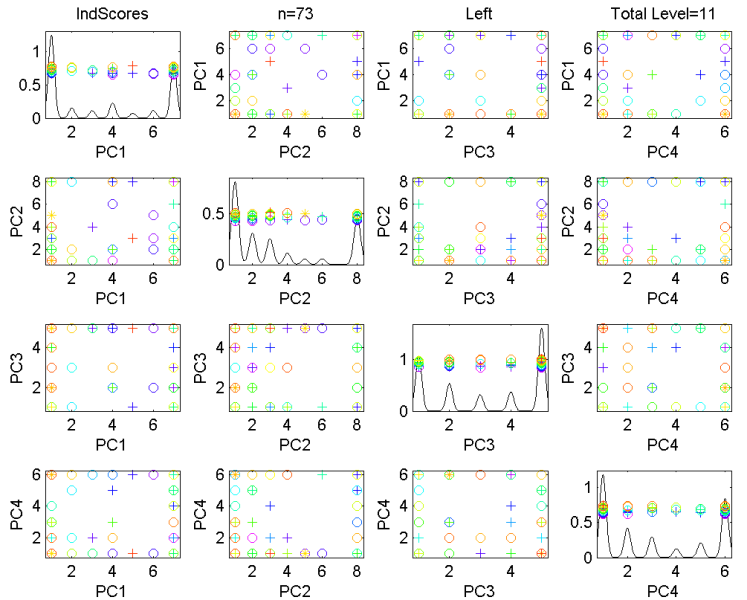}%
\caption{Scores Scatterplot for the Descendant Correspondence,
Left Side sub-population. Colors show age, symbols gender. No clear visual
patterns are apparent.}%
\label{Fig8}%
\end{center}
\end{figure}

Figure \ref{Fig8} shows the scores scatterplot for the set of left trees,
based on the descendant correspondence. The data points have been colored in
Figure \ref{Fig8}, to indicate age, which is an important covariate, as
discussed in Bullitt et al (2008). The color scheme
starts with purple for the youngest person (age $20$) and extends through a
rainbow type spectrum (blue-cyan-green-yellow-orange) to red for the oldest
(age $72$). An additional covariate, of possible interest, is sex, with
females shown as circles, males as plus signs, and two transgender cases
indicated using asterisks.

\begin{figure}
[h]
\begin{center}
\includegraphics[
natheight=6.916800in,
natwidth=10.666600in,
height=3.7766in,
width=5.1232in
]%
{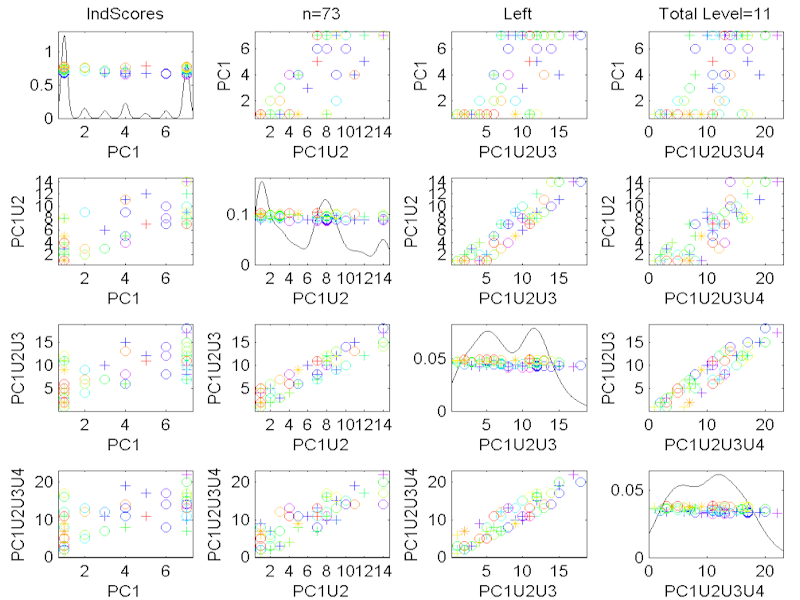}%
\caption{Cumulative Scores Scatterplot for the Descendant Correspondence,
Left Side sub-population.}
\label{Fig9}
\end{center}
\end{figure}

It was hoped that this visualization would reveal some interesting structure
with respect to age (color), but it is not easy to see any such connection in
Figure \ref{Fig8}. One reason for this is that the tree-lines only allow the
very limited range of scores as integers in the range $1$-$10$. A simple way to
generate a wider range of scores is to project not just onto simple tree-lines,
but instead onto their union, as defined in \eref{projunion}.
Figure \ref{Fig9}\ shows a scatterplot matrix, of several union PC scores, in
particular $PC1$ vs. $PC1\cup2$ (shorthand for $PC1\cup PC2$) vs.
$PC1\cup2\cup3$ vs. $PC1\cup2\cup3\cup4$. This combined plot, called the
\textit{cumulative scores scatterplot}, shows a better separation of the data
than is available in Figure \ref{Fig8}. The PC
unions show a banded structure, which again is an artifact that follows from
each PC score individually having a very limited range of possible values.
This seems to be a serious limitation of the tree-line approach to analyzing
population structure.

As with Figure \ref{Fig8}, there is unfortunately no readily apparent visual
connection between age and the visible population structure. However, visual
impression of this type can be tricky, and in particular it can be hard to see
some subtle effects.

\begin{figure}
[h]
\begin{center}
\includegraphics[
natheight=6.916800in,
natwidth=10.666600in,
height=3.7766in,
width=5.1232in
]
{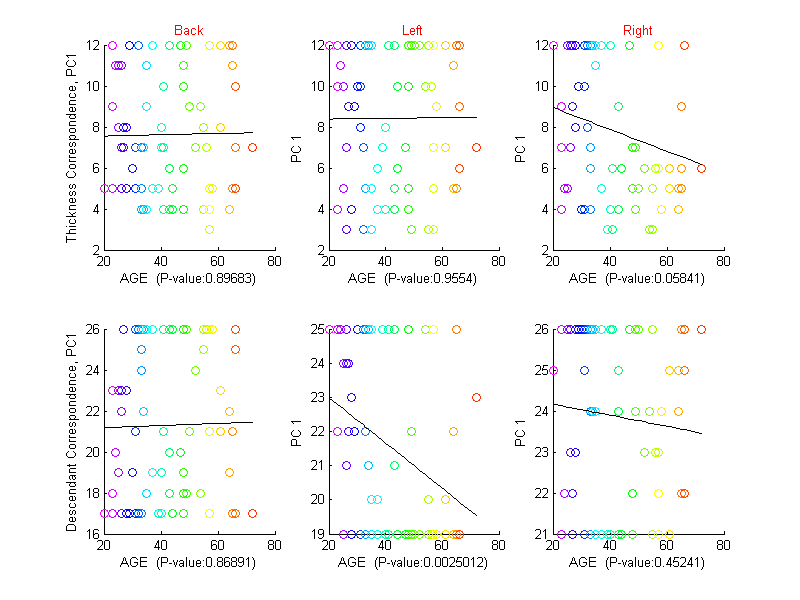}
\caption{Scatterplot of $PC1$ score versus age. Least squares fit regression
line suggests a downward trend in age. Trend is confirmed by the $p$-value of
$0.003$ (for significance of slope of the line).}
\label{Fig10}
\end{center}
\end{figure}

Figure \ref{Fig10} shows a view that more deeply scrutinizes the dependence of
the $PC1$ score on age, using a scatterplot, overlaid with the least squares
regression fit line. Note that most of the lines slope downwards, suggesting that
older people tend to have a smaller $PC1$ projection than younger people.
Statistical significance of this downward slope is tested by calculating the
standard linear regression $p$-value for the null hypothesis of $0$ slope. For the left tree, using the descendant
correspondence, the $p$-value is $0.0025$. This result is strongly
significant, indicating that this component is connected with age. This is
consistent with the results of Bullitt et al (2008), who noted a decreasing
trend with age in the total number of nodes. Our result is the first
location specific version of this.

Similar score versus age plots have been made, and hypothesis tests have been
run, for other PC components, and the resulting $p$-values, for the left tree using the descendent correspondence are summarized in
this table:

\[
\begin{array}
[c]{cc}
PC1 & 0.003\\
PC2 & 0.169\\
PC3 & 0.980\\
PC4 & 0.2984\\
PC1\cup2 & 0.003\\
PC1\cup2\cup3 & 0.004\\
PC1\cup2\cup3\cup4 & 0.007
\end{array}
\]

Note that for the individual PCs, only $PC1$ gives a statistically significant
result. For the cumulative PCs, all are significant,
but the significance diminishes as more components are added. This suggests
that it is really $PC1$ which is the driver of all of these results.

To interpret these results, recall from Figure \ref{Fig7}, that for the left
trees, $PC1$ chooses the left child for the first $3$ splits, and the right child
at the $4$th split. This suggests that there is not a significant difference
between the ages in the tree levels closer to the root, however, the difference does show
up when one looks at the deeper tree structure, in particular after the $4$th
split. This is consistent with the above remark, that for the left brain
sub-population, the first few splits did not seem to contain relevant
population information. Instead the effects of age only appear on splits
after level $4$.

We did a similar analysis of the back and right brain location
sub-populations, but none of these found significant results, so they are not
shown here. However, these can be found at the web site \cite{web}.

We also considered parallel results for the thickness correspondence, which
again did not yield significant results (but these are on the web site \cite{web}).
The fact that descendant correspondence gave some significant results, while
thickness never did, is one more indication that descendant correspondence is preferred.

\begin{figure}[h]
\begin{center}
\includegraphics[
natheight=9.375400in,
natwidth=12.500000in,
height=3.7766in,
width=5.028in
]
{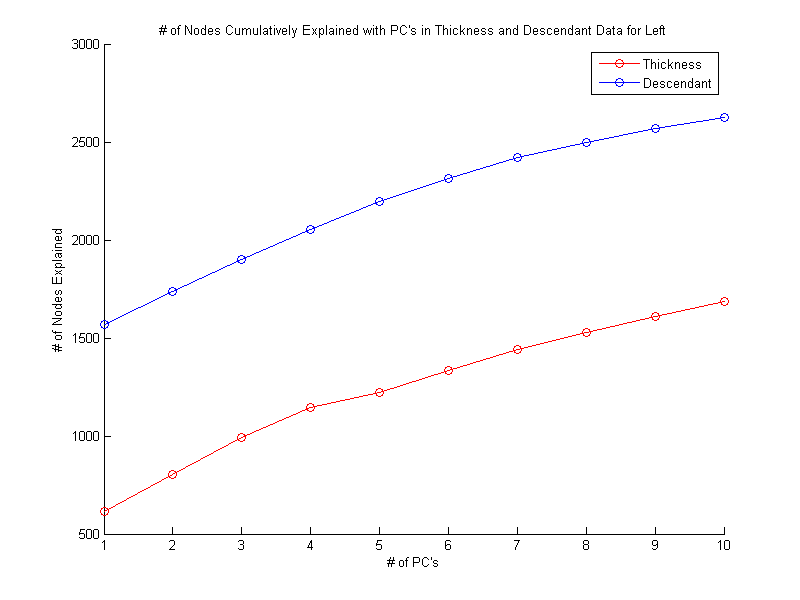}
\caption{Total number of nodes explained, as a function of Cumulative PC
Number. Shows that the descendant correspondence allows PCA to explain a
much higher proportion of the variation in the population than the thickness
correspondence.}
\label{Fig11}
\end{center}
\end{figure}

One more approach to the issue of correspondence choice is shown in Figure
\ref{Fig11}. This shows the amount of variation explained, as a function of
the order of the Cumulative Union PC, for both the thickness and the
descendant correspondences, for the left brain location sub-population. The
\textit{amount of variation explained} is defined to be the sum, over all
trees in the sub-population of the lengths of the projections. There are $5023$ nodes in total
for both correspondences. (The correspondence difference affects the locations of nodes, total count remains the same.)

It is not surprising that these curves are concave, since the first PC is designed to
explain the most variation, which each succeeding component explaining a
little bit less. But the important lesson from Figure \ref{Fig11} is that the
descendant correspondence allows PCA to explain much more
population structure, at each step, than the thickness correspondence.

In summary, there are several important consequences of this work:

\begin{itemize}
\item In real data sets with branching structure, tree PCA can reveal interesting insights, such as symmetry.

\item The descendant correspondence is clearly superior to the thickness
correspondence, and is recommended as the default choice in future studies.

\item As expected, the back sub-population is seen to have a more symmetric structure.

\item For the left sub-population there is a statistically significant structural age effect.

\item There seems to be room for improvement of the tree-line idea for doing
PCA on populations of trees. A possible improvement is to allow a richer branching structure, such as adding the next node as a child of one of the
last $2$ or $3$ nodes. We are exploring this methodology in our current research.
\end{itemize}

\section{Optimization proofs\label{Sec:Opt}}

This section is devoted to the proof of Theorem \ref{main} with some accompanying claims.

\begin{claim} \label{CLfirst} Let $L = \{ \ell_0, \dots, \ell_m \}$ be a tree-line, and $t$ a data tree. Then
\begin{eqnarray} \label{CLfirst1}
P_L(t)=\ell_0 \cup (t\cap V_L).
\end{eqnarray}
\end{claim}

\noindent\textbf{Proof:}\ \  Since
$\ell_i = \ell_{i-1} \cup v_i, \,$ we have
\begin{equation} \label{dtli}
\begin{array}{rcl}
d(t,\ell_i)
&=&\left\{\begin{array}{ll} d(t,\ell_{i-1})-1 &\ \mbox{if $v_i\in t$};\\
d(t,\ell_{i-1})+1 &\ \mbox{otherwise}.\end{array}\right. \\
\end{array}
\end{equation}
In other words, the distance of the tree to the line decreases as we keep adding nodes of $V_L$ that are in $t$, and
when we step out of $t$, the distance begins to increase, so Claim \eref{CLfirst1} follows.

\qed

\begin{claim} \label{ProjUnion} Let $L_1, \dots, L_q$ be tree-lines with a common starting point, and $t$ a data tree. Then
\begin{eqnarray*}
P_{L_1 \cup \dots \cup L_q}(t)=P_{L_1}(t) \cup \dots \cup P_{L_q}(t).
\end{eqnarray*}
\end{claim}

\noindent\textbf{Proof:} \ \ For simplicity, we only prove the statement for $q=2$. Assume that
\begin{eqnarray*}
    L_1 &=& \{ \ell_{1,0}, \ell_{1,1}, \dots, \ell_{1,p_1} \} \\
    L_2 &=& \{ \ell_{2,0}, \ell_{2,1}, \dots, \ell_{2,p_2} \}
\end{eqnarray*}
with $\ell_0 = \ell_{1,0} = \ell_{2,0},$  and
\beq
V_{L_1} = \{ v_{1,1}, \dots, v_{1,p_1} \}, \, V_{L_2} = \{ v_{2,1}, \dots, v_{2,p_2} \}.
\eeq
Also assume
\beqa \label{pl1t}
    P_{L_1}(t) &=& \ell_{1,r_1}, \\ \label{pl2t}
    P_{L_2}(t) &=& \ell_{2,r_2}. \;
\eeqa
For brevity, let us define
\beqa
f(i,j) & = & d(t, \ell_{1,i} \cup \ell_{2,j}) \; \text{for} \; 1 \leq i \leq p_1, \, 1 \leq j \leq p_2.
\eeqa
Using Claim \ref{CLfirst}, \eref{pl1t} means
\beqa \label{v1r}
v_{1,i} \in t, \; \text{if} \; i \leq r_1, \; \text{and} \; v_{1,i} \not\in t, \; \text{if} \; i > r_1,
\eeqa
hence
\beq \label{fij1}
\ba{rcl}
f(i,j) & \leq & f(i-1,j) \; \text{if} \; i \leq r_1; \;  \\
f(i,j) & \geq & f(i-1,j) \; \text{if} \; i > r_1.
\ea
\eeq
By symmetry, we have
\beq \label{fij2}
\ba{rcl}
f(i,j) & \leq & f(i,j-1) \; \text{if} \; j \leq r_2; \\
f(i,j) & \geq & f(i,j-1) \; \text{if} \; j > r_2.
\ea
\eeq
Overall, \eref{fij1} and \eref{fij2} imply that the function $f$ attains its minimum at $i=r_1, \, j=r_2, \,$ which is what we had to
prove.
\qed

\begin{claim} \label{findL} Let $S$ be a subset of $\Supp(T)$ which contains $\ell_0$.
For $v \in \Supp(T)$ define
\beq
\ba{rcl}
w_S(v) & = & \left\{ \ba{rl}  0, & \text{if} \,\, v \in S, \\
                              \sum_{v \in t_i} 1, & \text{otherwise}
                      \ea \right.
\ea
\eeq
Then among the treelines with starting tree $\ell_0$ the one which maximizes
$$
\sum_{t_i \in T} |(V_L \cup S) \cap t_i|
$$
is the one whose path $V_L$ maximizes the sum of the $w_S$ weights: $\sum_{v \in V_L} w_S(v)$.
\end{claim}
\noindent\textbf{Proof:} \ \ For $v \in \Supp(T), \,$ and a subtree $t$ of $\Supp(T), \,$ let us define
\beq
\ba{rcl}
\delta(v,t) & = & \left\{ \ba{rl}  1, & \text{if} \,\, v \in t, \\
                                   0, & \text{otherwise}
                      \ea \right.
\ea
\eeq
Then
\co{
\begin{eqnarray*}
\aml  \sum_{t_i \in T} |(V_L \cup S) \cap t_i| &=& \aml \sum_{t_i \in T} \sum_{v \in V_L \cup S} \delta(v, t_i) \\
                             & = &  \underset{L}{\arg\max} \sum_{v \in V_L \cup S} \sum_{t_i \in T} \delta(v, t_i) \\
                             & = & \underset{L}{\arg\max} \sum_{v \in V_L \cup S} w_\emptyset(v).
\end{eqnarray*}
}

$$
\ba{rcl}
\amal \sum_{t_i \in T} |(V_L \cup S) \cap t_i| &=& \amal \sum_{t_i \in T} \sum_{v \in V_L \cup S} \delta(v, t_i) \\
                             & = &  \amal \sum_{v \in V_L \cup S} \sum_{t_i \in T} \delta(v, t_i) \\
                             & = & \amal \sum_{v \in V_L \cup S} w_\emptyset(v) \\
                             & = & \amal \sum_{v \in  V_L} w_S(v).
\ea
$$

\qed

Finally, we prove our main result:

\noindent\textbf{Proof of Theorem \ref{main}:} \ \ For better intuition, we first give a proof when  $k=1. \,$ Using
Claim \ref{CLfirst} in Definition \ref{def-l1star}, we get
\beqast
L_1^* & = & \underset{L}{\arg\min} \sum_{t_i \in T} d( t_i, \ell_0 \cup (t_i\cap V_L)).
\eeqast
Since
$V_L$ is disjoint from $\ell_0, \,$
\beqast
L_1^* & = & \underset{L}{\arg\max} \sum_{t_i \in T} | V_L \cap t_i |,
\eeqast
the statement  follows from Claim \ref{findL} with $S = \emptyset$.

We now prove the statement for general $k.$
For an arbitrary data tree $t$, and tree-line $L$, we have
\beq \label{plk}
\ba{rcl}
P_{L_1^* \cup \dots \cup L_{k-1}^* \cup L}(t) & = & P_{L_1^*}(t) \cup \dots \cup P_{L_{k-1}^*}(t) \cup P_{L}(t) \\
                  & = & \ell_0 \cup ( V_{L_1^*} \cap t) \cup \dots \cup ( V_{L_{k-1}^*} \cap t) \cup ( V_{L} \cap t) \\
                  & = & \ell_0 \cup [ ( V_{L_1^*} \cup  \dots \cup V_{L_{k-1}^*} \cup V_L) \cap t ],
\ea
\eeq
with the first equation from Claim \ref{ProjUnion}, the second from Claim \ref{CLfirst}, and the third straightforward.

Combining \eref{plk} with \eref{lkstar} we get
\beq \label{lkstar-other}
L_k^*=\underset{L}{\arg\min} \sum_{t_i \in T} d(t_i, \ell_0 \cup [ ( ( V_{L_1^*} \cup  \dots \cup V_{L_{k-1}^*} \cup V_L) \cap t_i ]).
\eeq
Again, the paths of $L_1^*, \dots, L_{k-1}^*$ and $L$ are disjoint from $\ell_0, \,$ so \eref{lkstar-other} becomes
\beq \label{lkstar-other-2}
L_k^*=\underset{L}{\arg\max} \sum_{t_i \in T} | ( V_{L_1^*} \cup  \dots \cup V_{L_{k-1}^*} \cup V_L) \cap t_i |,
\eeq
so the statement follows from Claim \ref{findL} with $S = V_{L_1^*} \cup  \dots \cup V_{L_{k-1}^*}$.
\qed

\end{document}